\documentclass[aos,MSNbibl,seceqn,citesort,dvips]{arximspdf}

\doi{10.1214/11-AOS903}
\volume{39}
\issue{5}
\pubyear{2011}
\firstpage{2383}
\lastpage{2409}

\makeatletter

\newtheorem{theorem}{Theorem}
\newtheorem{corollary}{Corollary}

\newproclaim{definition}{Definition}

\newtheorem{lemma}{Lemma}
\newtheorem{proposition}{Proposition}

\makeatother

\begin{document}
\begin{frontmatter}

\title{On adaptive inference and confidence bands}
\runtitle{Adaptive inference}

\begin{aug}
\author[A]{\fnms{Marc} \snm{Hoffmann}\ead[label=e1]{marc.hoffmann@ensae.fr}}
\and
\author[B]{\fnms{Richard} \snm{Nickl}\corref{}\ead[label=e2]{r.nickl@statslab.cam.ac.uk}}
\runauthor{M. Hoffmann and R. Nickl}
\affiliation{ENSAE-CREST and University of Cambridge}
\address[A]{ENSAE-CREST and CNRS UMR 8050\\
Laboratoire d'Analyse\\
\quad et de Math\'{e}matiques Appliqu\'{e}es \\
3, avenue Pierre Larousse \\
92245 Malakoff Cedex\\
France\\
\printead{e1}}
\address[B]{Statistical Laboratory\\
Department of Pure Mathematics \\
\quad and Mathematical Statistics\\
University of Cambridge\\
CB3 0WB Cambridge\\
United Kingdom\\
\printead{e2}} %adresu isvedimo komanda gale!
\end{aug}

% HISTORY:
\received{\smonth{9} \syear{2010}}
\revised{\smonth{3} \syear{2011}}

% ABSTRACT
%
\begin{abstract}
The problem of existence of adaptive confidence bands for an unknown
density $f$ that belongs to a nested scale of H\"{o}lder classes over
$\mathbb R$ or $[0,1]$ is considered. Whereas honest adaptive inference
in this problem is impossible already for a pair of H\"{o}lder balls
$\Sigma(r), \Sigma(s), r \ne s$, of fixed radius, a nonparametric
distinguishability condition is introduced under which adaptive
confidence bands can be shown to exist. It is further shown that this
condition is necessary and sufficient for the existence of honest
asymptotic confidence bands, and that it is strictly weaker than
similar analytic conditions recently employed in Gin\'{e} and Nickl
[\textit{Ann. Statist.} \textbf{38} (2010) 1122--1170]. The exceptional
sets for which honest inference is not possible have vanishingly small
probability under natural priors on H\"{o}lder balls $\Sigma(s)$. If no
upper bound for the radius of the H\"{o}lder balls is known, a price
for adaptation has to be paid, and near-optimal adaptation is possible
for standard procedures. The implications of these findings for a
general theory of adaptive inference are discussed.
\end{abstract}

% KEYWORDS
%

\end{frontmatter}

%s1 ###
\section{Introduction}\label{sec1}

One of the intriguing problems in the paradigm of adaptive
nonparametric function estimation as developed in the last two decades
is what one could call the ``\textit{hiatus}'' between estimation and
inference, or, to be more precise, between the existence of adaptive
risk bounds and the nonexistence of adaptive confidence statements. In
a nutshell the typical situation in nonparametric statistics could be
described as follows: one is interested in a functional parameter $f$
that could belong either to $\Sigma$ or to $\Sigma'$, two sets that can
be distinguished by a certain ``structural property,'' such as
smoothness, with the possibility that $\Sigma\subset\Sigma'$. Based
on a sample whose distribution depends on $f$, one aims to find a
statistical procedure that adapts to the unknown structural property,
that is, that performs optimally without having to know whether $f \in
\Sigma$ or $f \in\Sigma'$. Now while such procedures can often be
proved to exist, the statistician cannot take advantage of this
optimality for inference: To cite Robins and van der Vaart \cite{RV06},
``An adaptive estimator can adapt to an underlying model, but does not
reveal which model it adapts to, with the consequence that
nonparametric confidence sets are necessarily much larger than the
actual discrepancy between an adaptive estimator and the true
parameter.''

We argue in this article that adaptive inference is possible if the
structural property that defines $\Sigma$ and $\Sigma'$ is
statistically identifiable, by which we shall mean here that the
nonparametric hypotheses $H_0\dvtx f \in\Sigma$ and $H_1\dvtx f \in
\Sigma'
\setminus\Sigma$ are asymptotically consistently distinguishable (in
the sense of Ingster \cite{I86,I93,I02}). In common adaptation
problems this will necessitate that certain unidentified parts of the
parameter space be removed, in other words, that the alternative
hypothesis $H_1$ be restricted to a subset $\tilde\Sigma$ of $\Sigma'
\setminus\Sigma$. One is in turn interested in choosing $\tilde
\Sigma
$ as large as possible, which amounts to imposing minimal
identifiability conditions on the parameter space. We shall make these
ideas rigorous in one key example of adaptive inference: confidence
bands for nonparametric density functions $f$ that adapt to the unknown
smoothness of $f$. The general approach, however, is not specific to
this example as we shall argue at the end of this introduction, and the
heuristic mentioned above is valid more generally.

The interest in the example of confidence bands comes partly from the
fact that the discrepancy between estimation and inference in this case
is particularly pronounced. Let us highlight the basic problem in a
simple ``toy adaptation'' problem. Consider $X_1,\ldots,X_n$ independent
and identically distributed random variables taking values in $[0,1]$
with common probability density function $f$ and joint law $\Pr_f$. We
are interested in the existence of confidence bands for $f$ that are
adaptive over two nested balls in the classical H\"{o}lder spaces
$\mathcal C^s([0,1]) \subset\mathcal C^r([0,1])$, $s>r$, of smooth
functions with norm given by \mbox{$\|\cdot\|_{s, \infty}$}; see Definition
\ref{besov} below. Define the class of densities
%
%e1.1 ###
%
\begin{equation} \label{sig}
\Sigma(s) := \Sigma(s, B) = \biggl\{f\dvtx[0,1] \to[0,\infty), \int_0^1
f(x)\,dx=1, \|f\|_{s,\infty} \le B \biggr\}\hspace*{-22pt}
\end{equation}
and note that $\Sigma(s) \subset\Sigma(r)$ for $s > r$. We shall
assume throughout that $B\ge1$ to ensure that $\Sigma(s)$ is nonempty.

A confidence band $C_n=C_n(X_1,\ldots,X_n)$ is a family of random
intervals
\[
\{C_n(y)=[c_n(y),c'_n(y)]\}_{y \in[0,1]}
\]
that contains graphs of densities $f\dvtx[0,1] \to[0, \infty)$. We denote
by $|C_n|={\sup_{y \in[0,1]}}|c'_n(y)-c_n(y)|$ the\vspace*{1pt} maximal diameter of
$C_n$. Following Li \cite{L89} the band $C_n$ is called asymptotically
\textit{honest} with level $\alpha$ for a family of probability
densities $\mathcal P$ if it satisfies the asymptotic coverage inequality
%
%e1.2 ###
%
\begin{equation} \label{hon}
\liminf_n \inf_{f \in\mathcal P} {\Pr}_f \bigl(f(y) \in C_n(y)\mbox{ }
\forall y \in[0,1] \bigr) \ge1- \alpha.
\end{equation}
We shall usually only write ${\Pr}_f(f \in C_n)$ for the coverage
probability if no confusion may arise. Note that $\mathcal P$ may (and
later typically will have to) depend on the sample size $n$. Suppose
the goal is to find a confidence band that is honest for the class
\[
\mathcal P^{\mathrm{all}} := \Sigma(s) \cup\Sigma(r) = \Sigma(r)
\]
and that is simultaneously \textit{adaptive} in the sense that the
expected diameter $E_{f}|C_n|$ of $C_n$ satisfies, for every $n$ (large enough),
%
%e1.3 ###
%
\begin{equation} \label{adapt0}
\sup_{f \in\Sigma(s)} E_f|C_n| \le L r_n(s),\qquad \sup_{f \in\Sigma
(r)} E_f|C_n| \le L r_n(r),
\end{equation}
where $L$ is a finite constant independent of $n$ and where
\[
r_n(s) = \biggl(\frac{\log n}{n} \biggr)^{{s}/({2s+1})}.
\]
Indeed even if $s$ were known no band could have expected diameter of
smaller order than $r_n(s)$ uniformly over $\Sigma(s)$ (e.g.,
Proposition \ref{nob} below), so that we are looking for a band that is
asymptotically honest for $\mathcal P^{\mathrm{all}}$ and that shrinks at the
fastest possible rate over $\Sigma(s)$ and $\Sigma(r)$ simultaneously.
It follows from Theorem~2 in Low \cite{L97} (see also
\cite{CL04,GW08}) that such bands do not exist.

\begin{theorem}[(Low)]
\label{low}
Any confidence band $C_n$ that is
honest over $\mathcal P^{\mathrm{all}}$ with level $\alpha<1$ necessarily
satisfies
\[
\lim_n \sup_{f \in\Sigma(s)}\frac{E_f|C_n|}{r_n(s)} = \infty.
\]
\end{theorem}

The puzzling fact is that this is in stark contrast to the situation in
estimation: adaptive estimators $\hat f_n$ such as those based on
Lepski's method \cite{L90} or wavelet thresholding \cite{DJKP96} can be
shown to satisfy simultaneously
\[
\sup_{f \in\Sigma(s)} E_f\|\hat f_n -f\|_\infty=O(r_n(s)),\qquad \sup
_{f \in\Sigma(r)} E_f\|\hat f_n -f\|_\infty=O(r_n(r));
\]
see \cite{GN09a,GN09b,GN10b} and Theorem \ref{supad} below. So while
$\hat f_n$ adapts to the unknown smoothness $s$, Theorem \ref{low}
reflects the fact that knowledge of the smoothness is still not
accessible for the statistician.

Should we therefore abstain from using adaptive estimators such as
$\hat f_n$ for inference? Gin\'{e} and Nickl \cite{GN10a} recently
suggested a new approach to this problem, partly inspired by Picard and
Tribouley \cite{PT00}. In \cite{GN10a} it was shown that one can
construct confidence bands $C_n$ and subsets $\bar\Sigma(\varepsilon,r)
\subset\Sigma(r)$, defined by a concrete analytical condition that
involves the constant $\varepsilon>0$, such that $C_n$ is asymptotically
honest for
\[
\mathcal P_\varepsilon= \Sigma(s) \cup\bar
\Sigma(\varepsilon,r)
\]
for every fixed $\varepsilon>0$, and such that $C_n$
is adaptive in the sense of (\ref{adapt0}). Moreover, these subsets
were shown to be topologically generic in the sense that the set
\[
\{f \in\Sigma(r) \mbox{ but } f \notin\bar\Sigma(\varepsilon,r)
\mbox{ for any } \varepsilon>0 \}
\]
that was removed is
\textit{nowhere dense} in the H\"{o}lder norm topology of $\mathcal
C^r$ (in fact in the relevant trace topology on densities). This says
that the functions $f \in\mathcal P^{\mathrm{all}}$ that prevent adaptation in
Theorem \ref{low} are in a certain sense negligible.

In this article we shall give a more statistical interpretation of
when, and if, why, adaptive inference is possible over certain subsets
of H\"{o}lder classes. Our approach will also shed new light on why
adaptation is possible over the sets $\bar\Sigma(\varepsilon, r)$.
Define, for $s>r$, the following class:
%
%e1.4 ###
%
\begin{equation} \label{supsep}\quad
\tilde\Sigma(r, \rho_n) := \tilde\Sigma(r, s,\rho_n, B) = \Bigl\{f
\in\Sigma(r,B)\dvtx\inf_{g \in\Sigma(s)}\|g-f\|_\infty\ge\rho_n
\Bigr\},
\end{equation}
where $\rho_n$ is a sequence of nonnegative real numbers. Clearly $\tilde
\Sigma(r, 0) = \Sigma(r)$, but if $\rho_n>0$, then we are removing
those elements from $\Sigma(r)$ that are not separated away from
$\Sigma
(s)$ in sup-norm distance by at least $\rho_n$. Inspection of the proof
of Theorem \ref{2class} shows that the set removed from $\Sigma(r)
\setminus\Sigma(s)$ is nonempty as soon as $\rho_n>0$.

Similar to above we are interested in finding a confidence band that is
honest over the class
\[
\mathcal P(\rho_n) := \Sigma(s) \cup\tilde\Sigma(r, \rho_n),
\]
and that is \textit{adaptive} in the sense of (\ref{adapt0}), in fact
only in the sense that
%
%e1.5 ###
%
\begin{equation} \label{adapt}
\sup_{f \in\Sigma(s)}E_f|C_n| \le L r_n(s), \qquad\sup_{f \in\tilde
\Sigma(r, \rho_n)}E_f|C_n| \le L r_n(r)
\end{equation}
for every $n$ (large enough). We know from Low's results that this is
impossible if $\rho_n=0$, but the question arises as to whether this
changes if $\rho_n>0$, and if so, what the smallest admissible choice
for $\rho_n$ is.

It was already noted or implicitly used in
\cite{HL02,JL03,B04,CL06,RV06} that there is a generic connection
between adaptive confidence
sets and minimax distinguishability of certain nonparametric
hypotheses. In our setting consider, for instance, testing the hypothesis
\[
H_0\dvtx f_0=1 \quad\mbox{against}\quad H_1\dvtx f_0 \in\mathcal
M,\qquad
\mathcal M \mbox{ finite}, \mathcal M \subset\tilde\Sigma(r, \rho_n).
\]
As we shall see in the proof of Theorem \ref{2class} below, an
adaptive confidence band over $\mathcal P(\rho_n)$ can be used to test
any such hypothesis consistently, and intuitively speaking an adaptive
confidence band should thus only exist if $\rho_n$ is of larger order
than the minimax rate of testing between $H_0$ and $H_1$ in the sense
of Ingster \cite{I86,I93}; see also the monograph \cite{I02}. For
confidence bands a natural separation metric is the supremum-norm (see,
however, also the discussion in the last paragraph of the
\hyperref[sec1]{Introduction}), and an exploration of the
corresponding testing problems
gives our main result, which confirms this intuition and shows moreover
that this lower bound is sharp up to constants at least in the case
where $B$ is known.

\begin{theorem} \label{2class} Let $s>r>0$. An adaptive and honest
confidence band over
\[
\Sigma(s) \cup\tilde\Sigma(r, \rho_n)
\]
exists if and only if $\rho_n$ is greater than or equal to the minimax
rate of testing between $H_0\dvtx f_0 \in\Sigma(s)$ and $H_1\dvtx f_0
\in
\tilde\Sigma(r, \rho_n)$, and this rate equals $r_n(r)$. More precisely:

\begin{longlist}[(a)]
\item[(a)] Suppose that $C_n$ is a confidence band that is asymptotically
honest with level $\alpha< 0.5$, over $\Sigma(s) \cup\tilde\Sigma(r,
\rho_n)$ and that is adaptive in the sense of (\ref{adapt}). Then
necessarily
\[
\liminf_n \frac{\rho_n}{r_n(r)}>0.
\]

\item[(b)] Suppose $B,r,s$ and $0<\alpha<1$ are given. Then there exists a
sequence $\rho_n$ satisfying
\[
\limsup_n \frac{\rho_n}{r_n(r)} < \infty
\]
and a confidence band $C_n=C_n(B,r,s,\alpha; X_1,\ldots, X_n)$ that is
asymptotically honest with level $\alpha$ and adaptive over $\Sigma(s)
\cup\tilde\Sigma(r, \rho_n)$ in the sense of (\ref{adapt}).

\item[(c)] Claims \textup{(a)} and \textup{(b)} still hold true if
$\Sigma(s)$ is replaced by the set
\[
\Bigl\{f \in\Sigma(s),\inf_{g \in\Sigma(t)} \|g-f\|_\infty\ge
Br_n(s)/2\Bigr\}
\]
for any $t>s$.
\end{longlist}
\end{theorem}

The last claim shows that the situation does not change if one removes
similar subsets from the smaller H\"{o}lder ball $\Sigma(s)$, in
particular removing the standard null-hypothesis $f_0=1$ used in the
nonparametric testing literature, or other very smooth densities,
cannot improve the lower bound for $\rho_n$.

Part (b) of Theorem \ref{2class} implies the following somewhat curious
corollary: since any $f \in\Sigma(r) \setminus\Sigma(s)$ satisfies
$\inf_{g \in\Sigma(s)}\|g-f\|_\infty>0$ (note that $\Sigma(s)$ is
\mbox{$\|
\cdot\|_\infty$}-compact), we conclude that $f \in\tilde\Sigma(r,
Lr_n(r))$ for every $L>0$, $n \ge n_0(f,r,L)$ large enough. We thus
have:

\begin{corollary} \label{dishon}
There exists a ``dishonest'' adaptive confidence band $C_n:=C_n(B, r,s,
\alpha; X_1,\ldots, X_n)$ that has asymptotic coverage for every fixed
$f \in\mathcal P^{\mathrm{all}}$; that is, $C_n$ satisfies
\[
\liminf_n {\Pr}_{f}(f \in C_n) \ge1-\alpha\qquad\forall f \in
\mathcal P^{\mathrm{all}}
\]
and
\begin{eqnarray*}
f &\in&\Sigma(s) \quad\Rightarrow\quad E_f|C_n| =O(r_n(s)),\\
f &\in&\Sigma(r)
\quad\Rightarrow\quad E_f|C_n| =O(r_n(r)).
\end{eqnarray*}
\end{corollary}

A comparison to Theorem \ref{low} highlights the subtle difference
between the minimax paradigm and asymptotic results that hold pointwise
in $f$: if one relaxes ``honesty,'' that is, if one removes the infimum
in (\ref{hon}), then Low's impossibility result completely disappears.
Note, however, that the index $n$ from which onwards coverage holds in
Corollary \ref{dishon} depends on $f$, so that the asymptotic result
cannot be confidently used for inference at a fixed sample size. This
is a reflection of the often neglected fact that asymptotic results
that are pointwise in $f$ have to be used with care for statistical
inference; see \cite{BLZ97,LP06} for related situations of this kind.

In contrast to the possibly misleading conclusion of Corollary \ref
{dishon}, Theorem \ref{2class} characterizes the boundaries of
``honest'' adaptive inference, and several questions arise.

\begin{longlist}
\item What is the relationship between the sets $\tilde
\Sigma(r, \rho_n)$ from Theorem \ref{2class} and the classes $\bar
\Sigma(\varepsilon,r)$ considered in \cite{GN10a}? Moreover, is
there a
``Bayesian'' interpretation of the exceptional sets that complements the
topological one?

\item The typical adaptation problem is not one over
two classes, but over a scale of classes indexed by a possibly
continuous smoothness parameter. Can one extend Theorem \ref{2class} to
such a setting and formulate natural, necessary and sufficient
conditions for the existence of confidence bands that adapt over a
continuous scale of H\"{o}lder classes?

\item Can one construct ``practical'' adaptive
nonparametric confidence bands? For instance, can one use bands that
are centered at wavelet or kernel estimators with data-driven
bandwidths? In particular can one circumvent having to know the radius
$B$ of the H\"{o}lder balls in the construction of the bands?
\end{longlist}

We shall give some answers to these questions in the remainder of the
article, and summarize our main findings here.

About question (i): we show in Proposition \ref{gnkill} that the
``statistical'' separation of $\Sigma(r)$ and $\Sigma(s)$ using the
sup-norm distance as in (\ref{supsep}) enforces a weaker condition on
$f \in\Sigma(r)$ than the analytic approach in \cite{GN10a}, so that
the present results are strictly more general for fixed smoothness
parameters $s$. We then move on to give a Bayesian interpretation of
the classes $\tilde\Sigma(r, \rho_n)$ and $\bar\Sigma(\varepsilon
, r)$:
we show in Proposition \ref{unif} that a natural Bayesian prior arising
from ``uniformly'' distributing suitably scaled wavelets on $\Sigma(r)$
concentrates on the classes $\tilde\Sigma(r, \rho_n)$ and $\bar
\Sigma
(\varepsilon, r)$ with overwhelming probability.

About question (ii): if the radius $B$ of the H\"{o}lder balls involved
is known, then one can combine a natural testing approach with recent
results in \cite{GN09a,GN09b,GN10b} to prove the existence of
adaptive nonparametric confidence bands over a scale of H\"{o}lder
classes indexed by a grid of smoothness parameters that grows dense in
any fixed interval $[r,R] \subset(0, \infty)$ as $n \to\infty$; see
Theorems~\ref{2classex},~\ref{discad}.

A full answer to question (iii) lies beyond the scope of this paper.
Some partial findings that seem of interest are the following: note
first that our results imply that the logarithmic penalties that
occurred in the diameters of the adaptive confidence bands in
\cite{GN10a} are not necessary if one knows the radius $B$. On the other
hand we show in Proposition \ref{nob} that if the radius $B$ is
unknown, then a certain price in the rate of convergence of the
confidence band cannot be
%circumvented---this also justifies the
%conventional wisdom, dating back to Bickel and Rosenblatt \cite{BR73},
%that ``undersmoothing'' is necessary in the construction of confidence
%bands.
circumvented, as $B$ cannot reliably be estimated without additional
assumptions on the model. This partly justifies the practice of
undersmoothing in the construction of confidence bands, dating back
to Bickel and Rosenblatt \cite{BR73}.
It leads us to argue that near-adaptive confidence bands that
can be used in practice, and that do not require the knowledge of $B$,
are more likely to follow from the classical adaptive techniques, like
Lepski's method applied to classical kernel or wavelet estimators,
rather than from the ``testing approach'' that we employ here to prove
existence of optimal procedures.

To conclude: the question as to whether adaptive methods should be used
for inference clearly remains a ``philosophical'' one, but we believe
that our results shed new light on the problem. That full adaptive
inference is not possible is a consequence of the fact that the typical
smoothness classes over which one wants to adapt, such as H\"{o}lder
balls, contain elements that are indistinguishable from a testing point
of view. On the other hand H\"{o}lder spaces are used by statisticians
to model regularity properties of unknown functions $f$, and it may
seem sensible to exclude functions whose regularity is not
statistically identifiable. Our main results give minimal
identifiability conditions of a certain kind that apply in this
particular case.

Our findings apply also more generally to the adaptation problem
discussed at the beginning of this introduction with two abstract
classes $\Sigma, \Sigma'$. We are primarily interested in confidence
statements that Cai and Low \cite{CL04} coin \textit{strongly adaptive}
(see Section 2.2 in their paper) and in our case this corresponds
precisely to requiring (\ref{hon}) and (\ref{adapt0}). If $\Sigma,
\Sigma'$ are convex, and if one is interested in a confidence interval
for a linear functional of the unknown parameter, Cai and Low show that
whether strong adaptation is possible or not is related to the
so-called ``inter-class modulus'' between $\Sigma, \Sigma'$, and their
results imply that in several relevant adaptation problems strongly
adaptive confidence statements are impossible. The
``separation-approach'' put forward in the present article (following
\cite{GN10a}) shows how strong adaptation can be rendered possible at
the expense of imposing statistical identifiability conditions on
$\Sigma, \Sigma'$, as follows: one first proves existence of a
risk-adaptive estimator $\hat f_n$ over $\Sigma, \Sigma'$ in some
relevant loss function. Subsequently one chooses a functional $\mathbb
F\dvtx\Sigma\times\Sigma' \to\mathbb[0,\infty)$, defines the
nonparametric model
\[
\mathcal P_n := \Sigma\cup\Bigl\{f \in\Sigma' \setminus\Sigma\dvtx
\inf_{g
\in\Sigma} \mathbb F(g,f) \ge\rho_n\Bigr\}
\]
and derives the minimax rate $\rho_n$ of testing $H_0\dvtx f \in
\Sigma$
against the generally nonconvex alternative $\{f \in\Sigma' \setminus
\Sigma\dvtx\inf_{g \in\Sigma} \mathbb F(g,f) \ge\rho_n\}$. Combining
consistent tests for these hypotheses with $\hat f_n$ allows for the
construction of confidence statements under sharp conditions on $\rho
_n$. A merit of this approach is that the resulting confidence
statements are naturally compatible with the statistical accuracy of
the adaptive estimator used in the first place. An important question
in this context, which is beyond the scope of the present paper, is the
optimal choice of the functional $\mathbb F$: for confidence bands it
seems natural to take $\mathbb F(f,g)=\|f-g\|_\infty$, but formalizing
this heuristic appears not to be straightforward. In more general
settings it may be less obvious to choose $\mathbb F$. These remain
interesting directions for future research.

%s2 ###
\section{\texorpdfstring{Proof of Theorem \protect\ref{2class}
and further results}{Proof of Theorem 2 and further results}} \label{main}

Let $X_1,\ldots,X_n$ be i.i.d. with probability density $f$ on $T$ which
we shall take to equal either $T=[0,1]$ or \mbox{$T=\mathbb R$}. We shall use
basic wavelet theory \cite{M92,HKPT,CDV93} freely throughout this
article, and we shall say that the wavelet basis is $S$-\textit
{regular} if the corresponding scaling functions $\phi_k$ and wavelets
$\psi_k$ are compactly supported and $S$-times continuously
differentiable on $T$. For instance, we can take Daubechies wavelets of
sufficiently large order $N=N(S)$ on $T=\mathbb R$ (see \cite{M92}) or
on $T=[0,1]$ (Section 4 in \cite{CDV93}).

We define H\"{o}lder spaces in terms of the moduli of the wavelet
coefficients of continuous functions. The wavelet basis consists of the
translated scaling functions $\phi_k$ and wavelets $\psi_{lk} =
2^{l/2}\psi_k(2^l(\cdot))$, where we add the boundary corrected scaling
functions and wavelets in case $T=[0,1]$. If $T=\mathbb R$ the indices
$k,l$ satisfy $l \in\mathbb N \cup\{0\}$, $k \in\mathbb Z$, but if
$T=[0,1]$ we require $l \ge J_0$ for some fixed integer $J_0=J_0(N)$
and then $k = 1,\ldots, 2^l$ for the $\psi_{lk}$'s, $k=1,\ldots
,N<\infty$
for the $\phi_k$'s. Note that $\psi_{lk}= 2^{l/2}\psi(2^l(\cdot
)-k)$ for
a fixed wavelet $\psi$ if either $T=\mathbb R$ or if $\psi_{lk}$ is
supported in the interior of $[0,1]$. Write shorthand $\alpha
_k(h)=\int
h \phi_k$, $\beta_{lk}(h)= \int h \psi_{lk}$.
\begin{definition} \label{besov}
Denote by $C(T)$ the space of bounded continuous real-valued functions
on $T$, and let $\phi_k$ and $\psi_k$ be $S$-regular Daubechies scaling
and wavelet functions, respectively. For $s < S$, the H\"{o}lder space
$\mathcal C^{s}(T)$ (\mbox{$=$}$\mathcal C^s$ when no confusion may arise) is
defined as the set of functions
\[
\Bigl\{ f \in C(T) \dvtx\|f\|_{s, \infty}\equiv\max\Bigl(\sup_{k}
|\alpha_k(f)|, \sup_{k, l} 2^{l(s+1/2)}|\beta_{lk}(f)|\Bigr)
<\infty\Bigr\}.
\]
\end{definition}

Define, moreover, for $s>0, B \ge1$, the class of densities
%
%e2.1 ###
%
\begin{equation} \label{sig}\qquad
\Sigma(s) := \Sigma(s,B, T) = \biggl\{f\dvtx T \to[0,\infty), \int_{T}
f(x)\,dx=1, \|f\|_{s,\infty} \le B \biggr\}.
\end{equation}
It is a standard result in wavelet theory (Chapter 6.4 in \cite{M92}
for $T=\mathbb R$ and Theorem~4.4 in \cite{CDV93} for $T=[0,1]$) that
$\mathcal C^s$ is equal, with equivalent norms, to the classical H\"
{o}lder--Zygmund spaces $C^s$. For $T=\mathbb R$, $0<s<1$, these spaces
consist of all functions $f \in C(\mathbb R)$ for which $\|f\|_\infty
+\sup_{x \ne y, x,y \in\mathbb R} (|f(x)-f(y)|/|x-y|^s)$ is finite.
For noninteger $s>1$ the space $C^s$ is defined by requiring $D^{[s]}f$
of $f \in C(\mathbb R)$ to exist and to be contained in $C^{s-[s]}$.
The Zygmund class $C^1$ is defined by requiring
$|f(x+y)+f(x-y)-2f(x)|\le C|y|$ for all $x,y \in\mathbb R$, some
$0<C<\infty$ and $f \in C(\mathbb R)$, and the case $m<s\le m+1$
follows by requiring the same condition on the $m$th derivative of $f$.
The definitions for $T=[0,1]$ are similar; we refer to \cite{CDV93}.

Define the projection kernel $K(x,y)=\sum_k \phi_k(x)\phi_k(y)$ and write
\begin{eqnarray*}
K_j(f)(x) &=& 2^j \int_T K(2^jx,2^jy)f(y)\,dy \\
&=& \sum_k \alpha_k(f)\phi_k + \sum_{l=J_0}^{j-1} \sum_k \beta
_{lk}(f)\psi_{lk}
\end{eqnarray*}
for the partial sum of the wavelet series of a function $f$ at
resolution level $j\ge J_0+1$, with the convention that $J_0=0$ if
$T=\mathbb R$.

If $X_1,\ldots,X_n$ are i.i.d. $\sim f$ then an unbiased estimate of
$K_j(f)$ is, for $\hat\alpha_k = (1/n)\sum_{i=1}^n \phi_k(X_i),
\hat
\beta_{lk}= (1/n) \sum_{i=1}^n \psi_{lk}(X_i)$ the empirical wavelet
coefficients,
%
%e2.2 ###
%
\begin{equation} \label{est}
f_n(x,j)= \frac{2^j}{n} \sum_{i=1}^n K(2^jx, 2^jX_i) = \sum_k \hat
\alpha_k\phi_k + \sum_{l=J_0}^{j-1} \sum_k \hat\beta_{lk}\psi_{lk}.
\end{equation}

%s2.1 ###
\subsection{\texorpdfstring{Proof of Theorem \protect\ref{2class}}{Proof of Theorem 2}}

We shall first prove Theorem \ref{2class} to lay out the main ideas. We
shall prove claims (a) and (b), that this also solves the testing
problem $H_0\dvtx f_0 \in\Sigma(s)$ against $H_1\dvtx f_0 \in\tilde
\Sigma(r,
\rho_n)$ follows from the proofs. The proof of claim (c) is postponed
to Section \ref{rem}. Let us assume $B \ge2$ to simplify some
notation. Take $j_n^* \in\mathbb N$ such that
\[
2^{j_n^*} \simeq\biggl(\frac{n}{\log n}\biggr)^{{1}/({2r+1})}
\]
is satisfied, where $\simeq$ denotes two-sided inequalities up to
universal constants.

($\Leftarrow$): Let us show that $\liminf_n (\rho_n/r_n(r))=0$ leads to
a contradiction. In this case $\rho_n/r_n(r) \to0$ along a subsequence
of $n$, and we shall still index this subsequence by $n$. Let $f_0=1$
on $[0,1]$ and define, for $\varepsilon>0$, the functions
\[
f_m := f_0 + \varepsilon2^{-j(r+1/2)} \psi_{jm},
\]
where $m=1,\ldots,M, c_02^{j} \le M <2^j$, $j\ge0, c_0 >0$, and where
$\psi$ is a Daubechies wavelet of regularity greater than $s$, chosen
in such a way that $\psi_{jm}$ is supported in the interior of $[0,1]$
for every $m$ and $j$ large enough. (This is possible using the
construction in Theorem 4.4 in \cite{CDV93}.) Since $\int_0^1 \psi=0$
we have $\int_0^1 f_m=1$ for every $m$ and also $f_m \ge0$ $\forall m$
if $\varepsilon>0$ is chosen small enough depending only on~$\|\psi\|
_\infty$. Moreover, for any $t>0$, using the definition of \mbox{$\|\cdot\|
_{t, \infty}$} and since $c(\phi)\equiv\sup_k|\int_0^1 \phi_k| \le
\sup
_k\|\phi_k\|_2 = 1$,
%
%e2.3 ###
%
\begin{equation} \label{norm1}
\|f_m\|_{t,\infty} = \max\bigl(c(\phi),
\varepsilon2^{j(t-r)}\bigr),\qquad
m=1,\ldots,M,
\end{equation}
so $f_m \in\Sigma(r)$ for $\varepsilon\le2$ (recall $B \ge2$) and
every $j$ but $f_m \notin\Sigma(s)$ for $j$ large enough depending
only on $s,r, B, \varepsilon$.

Note next that
\[
|\beta_{lk} (h)| = \biggl| 2^{l/2} \int\psi_k(2^lx)h(x)\,dx \biggr|
\le 2^{-l/2}\|\psi_k\|_1 \|h\|_\infty\le2^{-l/2} \|h\|_\infty
\]
for every $l, k$, and any bounded function $h$ implies
%
%e2.4 ###
%
\begin{equation}\label{lbdu}
\|h\|_\infty\ge\sup_{l \ge0, k} 2^{l/2}|\beta_{lk}(h)|
\end{equation}
so that, for $g \in\Sigma(s)$ arbitrary,
%
%e2.5 ###
%
\begin{eqnarray} \label{suplbd-} \|f_m-g\|_\infty&\ge& \sup_{l \ge0,
k} 2^{l/2} |\beta_{lk}(f_m)-\beta_{lk}(g)| \nonumber\\
& \ge& \varepsilon2^{-jr} - 2^{j/2}|\beta_{jk}(g)| \ge\varepsilon2^{-jr}
- B2^{-js} \\
&\ge& \frac{\varepsilon}{2} 2^{-jr}\nonumber
\end{eqnarray}
for every $m$ and for $j \ge j_0, j_0=j_0(s,r, B, \varepsilon)$.
Summarizing we see that
\[
f_m \in\tilde\Sigma\biggl(r, \frac{\varepsilon}{2} 2^{-jr}\biggr)
\qquad\forall m=1,\ldots,M
\]
for every $j \ge j_0$. Since $\rho_n=o(r_n(r))$, $r_n(r) \simeq
2^{-j_n^*r}$, we can find $j_n>j_n^*$ such that
%
%e2.6 ###
%
\begin{equation} \label{rho-}
\rho'_n := \max(\rho_n, r_n(s) \log n) \le\frac{\varepsilon}{2} 2^{-j_nr}
=o(2^{-j_n^*r})
\end{equation}
in particular $f_m \in\tilde\Sigma(r,\rho_n')$ for every
$m=1,\ldots
,M$ and every $n \ge n_0, n_0=n_0(s,r,B, \varepsilon)$.

Suppose now $C_n$ is a confidence band that is adaptive and honest over
$\Sigma(s) \cup\tilde\Sigma(r, \rho_n)$, and consider testing
\[
H_0\dvtx
f = f_0 \quad\mbox{against}\quad H_1\dvtx f \in\{f_1,\ldots,
f_M\}=:\mathcal M.
\]
Define a test $\Psi_n$ as follows: if no $f_m \in
C_n$, then $\Psi_n=0$, but as soon as one of the $f_m$'s is contained in
$C_n$, then $\Psi_n=1$. We control the error probabilities of this
test. Using (\ref{suplbd-}), Markov's inequality, adaptivity of the
band, (\ref{rho-}) and noting $r_n(s) = o(\rho_n')$, we deduce
\begin{eqnarray*}
{\Pr}_{f_0}(\Psi_n \ne0) &=& {\Pr}_{f_0} (f_m \in C_n \mbox{ for some
} m)\\
&=& {\Pr}_{f_0} (f_m, f_0 \in C_n\mbox{ for some } m) \\
&&{} + {\Pr}_{f_0} (f_m \in C_n \mbox{ for some } m, f_0 \notin C_n) \\
&\le& {\Pr}_{f_0} (\|f_m-f_0\|_\infty\le|C_n| \mbox{ for some } m) +
\alpha+ o(1) \\
&\le& {\Pr}_{f_0}(|C_n| \ge\rho'_n) + \alpha+o(1)
\\
&\le& E_{f_0}|C_n|/\rho'_n + \alpha+o(1) = \alpha+ o(1).
\end{eqnarray*}
Under any alternative $f_m \in\tilde\Sigma(r, \rho_n')$, invoking
honesty of the band we have
\[
P_{f_m}(\Psi_n = 0) = {\Pr}_{f_m}(\mbox{no } f_k \in C_n) \le{\Pr
}_{f_m} (f_m \notin C_n) \le\alpha+ o(1)
\]
so that summarizing we have
%
%e2.7 ###
%
\begin{equation} \label{bul}
\limsup_n\Bigl(E_{f_0}\Psi_n + \sup_{f \in\mathcal M}E_{f}(1-\Psi
_n)\Bigr) \le2 \alpha<1.
\end{equation}
On the other hand, if $\tilde\Psi$ is \textit{any} test (any
measurable function of the sample taking values $0$ or $1$), we shall
now prove
%
%e2.8 ###
%
\begin{equation} \label{contr}
\liminf_n \inf_{\tilde\Psi} \Bigl(E_{f_0}\tilde\Psi+ \sup_{f \in
\mathcal M}E_{f}(1-\tilde\Psi)\Bigr) \ge1,
\end{equation}
which contradicts (\ref{bul}) and completes this direction of the
proof. The proof follows ideas in \cite{I86}. We have, for every $\eta>0$,
\begin{eqnarray*}
E_{f_0}\tilde\Psi+ \sup_{f \in\mathcal M} E_f(1-\tilde\Psi) &\ge&
E_{f_0}(1\{\tilde\Psi=1\}) + \frac{1}{M} \sum_{m=1}^M
E_{f_m}(1-\tilde
\Psi) \\
&\ge& E_{f_0}(1\{\tilde\Psi=1\} + 1\{\tilde\Psi=0\}Z ) \\
&\ge& (1-\eta) {\Pr}_{f_0}(Z \ge1- \eta),
\end{eqnarray*}
where $Z= M^{-1} \sum_{m=1}^M (dP^n_m/dP^n_0)$ with $P^n_m$ the product
probability measures induced by a sample of size $n$ from the density
$f_m$. By Markov's inequality,
\[
{\Pr}_{f_0}(Z \ge1- \eta) \ge1- \frac{E_{f_0}|Z-1|}{\eta} \ge1 -
\frac{\sqrt{E_{f_0}(Z-1)^2}}{\eta}
\]
for every $\eta>0$, and we show that the last term converges to zero.
Writing (in abuse of notation) $\gamma_j = \varepsilon2^{-j_n(r+1/2)}$,
using independence, orthonormality of $\psi_{jm}$ and $\int\psi
_{jm}=0$ repeatedly as well as $(1+x) \le e^x$, we see
\begin{eqnarray*}
E_{f_0}(Z-1)^2 &=& \frac{1}{M^2} \int_{[0,1]^n} \Biggl(\sum_{m=1}^M \Biggl(\prod
_{i=1}^n f_m(x_i) -1\Biggr) \Biggr)^2 \,dx \\
&=& \frac{1}{M^2} \int_{[0,1]^n} \Biggl(\sum_{m=1}^M \Biggl( \prod_{i=1}^n
\bigl(1+\gamma_j \psi_{jm}(x_i)\bigr) -1\Biggr) \Biggr)^2 \,dx \\
&=& \frac{1}{M^2} \sum_{m=1}^M \int_{[0,1]^n} \Biggl( \prod_{i=1}^n
\bigl(1+\gamma_j \psi_{jm}(x_i)\bigr) -1 \Biggr)^2 \,dx \\
&=&\frac{1}{M^2} \sum_{m=1}^M \Biggl( \int_{[0,1]^n} \prod_{i=1}^n
\bigl(1+\gamma_j \psi_{jm}(x_i)\bigr)^2 \,dx -1 \Biggr) \\
&=& \frac{1}{M^2} \sum_{m=1}^M \Biggl( \biggl( \int_{[0,1]} \bigl(1+\gamma_j
\psi_{jm}(x)\bigr)^2 \,dx \biggr)^n -1 \Biggr) \\
&=& \frac{1}{M} \bigl((1+\gamma_j^2)^n -1 \bigr) \le\frac{e^{n \gamma
_j^2}-1}{M}.
\end{eqnarray*}
Now using (\ref{rho-}) we see $n \gamma_j^2 = \varepsilon^2
n2^{-j_n(2r+1)} = o(\log n)$ so that $e^{n \gamma_j^2} = o(n^{\kappa})$
for every $\kappa>0$, whereas $M \simeq2^{j_n} \ge2^{j_n^*} \simeq
r_n(r)^{-1/r}$ still diverges at a fixed polynomial rate in $n$, so
that the last quantity converges to zero, which proves (\ref{contr})
since $\eta$ was arbitrary.

($\Rightarrow$): Let us now show that an adaptive band $C_n$ can be
constructed if $\rho_n$ equals $r_n(r)$ times a large enough constant,
and if the radius $B$ is known. The remarks after Definition \ref
{besov} imply that $\|f\|_\infty\le k\|f\|_{s, \infty} \le kB$ for
some $k>0$. Set
%
%e2.9 ###
%
\begin{equation} \label{sigma}
\sigma(j) := \sigma(n,j) := \sqrt{kB \frac{2^jj}{n}},\qquad \rho_n:= L'
\sigma(j^*_n) \simeq r_n(r)
\end{equation}
for $L'$ a constant to be chosen later. Using Definition \ref{besov}
and $\sup_x \sum_k|\psi_k(x)|<\infty$, we have for $f_n$ from (\ref
{est}) based on wavelets of regularity $S>s$
%
%e2.10 ###
%
\begin{equation} \label{bias}
\|E_ff_n(j_n^*)-f\|_\infty= \|K_{j_n^*}(f)-f\|_\infty\le b_0
2^{-j_n^*r} \le b \sigma(j^*_n)
\end{equation}
for some constants $b_0, b$ that depend only on $B, \psi$.

Define the test statistic $\hat d_n:= \inf_{g \in\Sigma(s)} \|
f_n(j^*_n)-g\|_\infty$.
Let now $\hat f_n(y)$ be any estimator for $f$ that is exact rate
adaptive over $\Sigma(s) \cup\Sigma(r)$ in sup-norm risk; that is,
$\hat f_n$ satisfies simultaneously, for some fixed constant $D$
depending only on $B,s,r$
%
%e2.11 ###
%
\begin{equation} \label{adap}\qquad
\sup_{f \in\Sigma(r)}E_f \|\hat f_n - f\|_\infty\le D r_n(r),\qquad
\sup_{f \in\Sigma(s)}E_f\|\hat f_n - f\|_\infty\le D r_n(s).
\end{equation}
Such estimators exist; see Theorem \ref{supad} below. Define the
confidence band $C_n \equiv\{C_n(y), y \in\mathbb[0,1]\}$ to equal
\[
\hat f_n(y) \pm Lr_n(r)\qquad\mbox{if }\hat d_n > \tau
\quad\mbox{and}\quad\hat f_n(y) \pm Lr_n(s)
\qquad\mbox{if }\hat d_n \le\tau,y \in[0,1],
\]
where $\tau= \kappa\sigma(j_n^*)$, and where $\kappa$ and $L$ are
constants to be chosen below.

We first prove that $C_n$ is an honest confidence band for $f \in
\Sigma
(s) \cup\tilde\Sigma(r, \rho_n)$ when $\rho_n$ is as above with $L'$
large enough depending only on $\kappa, B$. If $f \in\Sigma(s)$ we
have coverage since adaptivity of $\hat f_n$ implies, by Markov's inequality,
\begin{eqnarray*}
\inf_{f \in\Sigma(s)}{\Pr}_f (f \in C_n )&\ge& 1- \sup_{f
\in\Sigma(s)}{\Pr}_f \bigl( \|\hat f_n-f\|_\infty> Lr_n(s) \bigr) \\
& \ge& 1-\frac{1}{Lr_n(s)} \sup_{f \in\Sigma(s)}E_f \|\hat f_n -f\|
_\infty\\
&\ge& 1-\frac{D}{L},
\end{eqnarray*}
which\vspace*{1pt} can be made greater than $1-\alpha$ for any $\alpha>0$ by
choosing $L$ large enough depending only on $K, B, \alpha, r,s$. When
$f \in\tilde\Sigma(r, \rho_n)$ there is the danger of $\hat d_n \le
\tau$ in which case the size of the band is too small. In this case,
however, we have, using again Markov's inequality,
\[
\inf_{f \in\tilde\Sigma(r, \rho_n)}{\Pr}_f (f \in C_n) \ge
1 - \frac{\sup_{f \in\tilde\Sigma(r, \rho_n)}E_f\|\hat f_n -f\|
_\infty
}{Lr_n(r)} - \sup_{f \in\tilde\Sigma(r, \rho_n)}{\Pr}_f(\hat d_n
\le
\tau)
\]
and the first term subtracted can be made smaller than $\alpha$ for $L$
large enough in view of (\ref{adap}). For the second note that ${\Pr}_f
(\hat d_n \le\tau)$ equals, for every $f \in\tilde\Sigma(r, \rho_n)$,
\begin{eqnarray*}
&& {\Pr}_f \Bigl(\inf_{g \in\Sigma(s)} \|f_n(j_n^*)-g\|_\infty\le
\kappa\sigma(j_n^*) \Bigr) \\
&&\qquad \le{\Pr}_f \Bigl(\inf_g \|f-g\|_\infty- \|
f_n(j_n^*)-E_ff_n(j_n^*)\|_\infty\\
&&\qquad\quad\hspace*{70pt}{} - \|K_{j_n^*}(f)-f\|_\infty\le
\kappa\sigma(j_n^*) \Bigr) \\
&&\qquad \le{\Pr}_f \bigl(\rho_n - \|K_{j_n^*}(f)-f\|_\infty- \kappa\sigma
(j_n^*) \le\|f_n(j_n^*)-E_ff_n(j_n^*) \|_\infty\bigr) \\
&&\qquad \le{\Pr}_f \bigl(\|f_n(j_n^*)-E_ff_n(j_n^*) \|_\infty\ge
(L'-\kappa-b) \sigma(j_n^*) \bigr) \\
&&\qquad\le ce^{-cj_n^*} = o(1)
\end{eqnarray*}
for some $c>0$, by choosing $L'=L'(\kappa, B, K)$ large enough
independent of $f \in\tilde\Sigma(r, \rho_n)$, in view of
Proposition \ref{ineq} below. This completes the proof of coverage of
the band.

We now turn to adaptivity of the band and verify (\ref{adapt}). By
definition of $C_n$ we have almost surely
\[
|C_n| \le Lr_n (r),
\]
so the case $f \in\tilde\Sigma(r, \rho_n)$ is proved. If $f \in
\Sigma
(s)$ then, using (\ref{bias}) and Proposition~\ref{ineq},
\begin{eqnarray*}
E_f|C_n| &\le& L r_n(r) {\Pr}_f (\hat d_n > \tau) + L r_n(s) \\
&\le& Lr_n(r){\Pr}_f \Bigl(\inf_{g \in\Sigma(s)} \|f_n(j_n^*)-g\|
_\infty> \kappa\sigma(j_n^*) \Bigr) + Lr_n(s) \\
&\le& Lr_n(r){\Pr}_f \bigl(\|f_n(j_n^*)-f\|_\infty> \kappa\sigma
(j_n^*) \bigr) + Lr_n(s) \\
&\le& Lr_n(r) {\Pr}_f \bigl(\|f_n(j_n^*)-E_ff_n(j_n^*)\|_\infty>
(\kappa-b) \sigma(j_n^*) \bigr) + L r_n(s)\\
&\le& Lr_n(r) ce^{-cj_n^*} + Lr_n(s) = O(r_n(s))
\end{eqnarray*}
since $c$ can be taken sufficiently large by choosing $\kappa=\kappa
(K,B)$ large enough. This completes the proof of the second claim of
Theorem \ref{2class}.

%s2.2 ###
\subsection{Unknown radius $B$}

The existence results in the previous section are not entirely
satisfactory in that the bands constructed to prove existence of
adaptive procedures cannot be easily implemented. Particularly the
requirement that the radius $B$ of the H\"{o}lder ball be known is
restrictive. A first question is whether exact rate-adaptive bands
exist if $B$ is unknown, and the answer turns out to be no. This in
fact is not specific to the adaptive situation, and occurs already for
a fixed H\"{o}lder ball, as the optimal size of a confidence band
depends on the radius~$B$. The following proposition is a simple
consequence of the formula for the exact asymptotic minimax constant
for density estimation in sup-norm loss as derived in~\cite{KN99}.
\begin{proposition} \label{nob}
Let $X_1,\ldots,X_n$ be i.i.d. random variables taking values in $[0,1]$
with density $f \in\Sigma(r,B, [0,1])$ where $0<r<1$. Let $C_n$ be a
confidence band that is asymptotically honest with level $\alpha$ for
$\Sigma(r,B, [0,1])$. Then
\[
\liminf_n\sup_{f \in\Sigma(r,B, [0,1])} \frac{E_f|C_n|}{r_n(r)}
\ge
cB^p (1-\alpha)
\]
for some fixed constants $c,p>0$ that depend only on $r$.
\end{proposition}

In particular if $C_n$ does not depend on $B$, then $E_f|C_n|$ cannot
be of order $r_n(r)$ uniformly over $\Sigma(r,B,[0,1])$ for every
$B>0$, unless $B$ can be reliably estimated, which for the full H\"{o}lder
ball is impossible without additional assumptions.
%. The best we can therefore expect of a confidence band $C_n$ that
%does not require the knowledge of $B$ is that one has to pay an
%arbitrarily slowly divergent penalty in the rate of convergence of
%$E_f|C_n|$.
It can be viewed as one explanation for why undersmoothing
is necessary to construct ``practical'' asymptotic confidence bands.

%s2.3 ###
\subsection{Confidence bands for adaptive estimators} \label{ansep}

The usual risk-adaptive estimators such as those based on Lepski's
\cite{L90} method or wavelet thresholding \cite{DJKP96} do not
require the
knowledge of the H\"{o}lder radius $B$. As shown in \cite{GN10a} (see
also~\cite{KNP10}) such estimators can be used in the construction of
(near-)adaptive confidence bands under certain analytic conditions on
the elements of $\Sigma(s)$. Let us briefly describe the results in
\cite{GN10a,KNP10}. Let $\ell_{n}$ be a sequence of positive integers
(typically $\ell_n \to\infty$ as $n \to\infty$) and define, for $K$
the wavelet projection kernel associated to some $S$-regular wavelet
basis, $S>s$
%
%e2.12 ###
%
\begin{equation} \label{barsig}
\bar\Sigma(\varepsilon,s, \ell_n) := \{f \in\Sigma(s)\dvtx
\varepsilon
2^{-ls} \le\|K_l(f)-f\|_\infty\le B2^{-ls}\mbox{ } \forall l \ge\ell_n
\}.\hspace*{-22pt}
\end{equation}
The conditions in \cite{GN10a,KNP10} are slightly weaker in that they
have to hold only for $l \in[\ell_n, \ell_n']$ where $\ell'_n -
\ell_n
\to\infty$. This turns out to be immaterial in what follows, however,
so we work with these sets to simplify the exposition.

Whereas the upper bound in (\ref{barsig}) is automatic for functions in
$\Sigma(s)$, the lower bound is not. However one can show that a lower
bound on $\|K_l(f)-f\|_\infty$ of order $2^{-ls}$ is ``topologically''
generic in the H\"{o}lder space $\mathcal C^s(T)$. The following is
Proposition 4 in \cite{GN10a}.
\begin{proposition} \label{5} Let $K$ be $S$-regular with $S>s$. The
set
\[
\bigl\{f\mbox{: there exists no } \varepsilon>0, l_0 \ge0 \mbox{
s.t. } \|K_l(f)-f\|_\infty\ge\varepsilon2^{-l(s+1/2)}\mbox{ }\forall l \ge
l_0 \bigr\}
\]
is nowhere dense in the norm topology of $\mathcal C^s(\mathbb R)$.
\end{proposition}

Using this condition, \cite{GN10a} constructed an estimator $\hat f_n$
based on Lepski's method applied to a kernel or wavelet density
estimator such that
%
%e2.13 ###
%
\begin{equation} \label{lim0}
\hat A_n \biggl(\sup_{y \in[0,1]}\biggl|\frac{\hat f_n(y)-f(y)}{\hat
\sigma_n\sqrt{\hat f_n(y)}}\biggr|-\hat B_n \biggr) \to^d Z
\end{equation}
as $n \to\infty$, where $Z$ is a standard Gumbel random variable and
where $\hat A_n, \hat B_n, \hat\sigma_n$ are some random constants. If
$\ell_n$ is chosen such that
%
%e2.14 ###
%
\begin{equation} \label{R}
2^{\ell_n} \simeq\biggl(\frac{n}{\log n}\biggr)^{1/(2R+1)},
\end{equation}
then the limit theorem (\ref{lim0}) is uniform in relevant unions over
$s \in[r,R], r>0$, of H\"{o}lder classes $\bar\Sigma(\varepsilon,
s, \ell
_n)$. Since the constants $\hat A_n, \hat B_n, \hat\sigma_n$ in (\ref
{lim0}) are known, confidence bands can be retrieved directly from the
limit distribution, and \cite{GN10a} further showed that so-constructed
bands are near-adaptive: they shrink at rate $O_P(r_n(s) u_n)$ whenever
$f \in\bar\Sigma(\varepsilon,s, \ell_n)$, where $u_n$ can be taken of
the size $\log n$. See Theorem 1 in \cite{GN10a} for detailed
statements. As shown in Theorem 4 in \cite{KNP10}, the restriction $u_n
\simeq\log n$ can be relaxed to $u_n \to\infty$ as $n \to\infty$, at
least if one is not after exact limiting distributions but only after
asymptotic coverage inequalities, and this matches Proposition \ref
{nob}, so that these bands shrink at the optimal rate in the case where
$B$ is unknown.

Obviously it is interesting to ask how the sets in (\ref{barsig})
constructed from analytic conditions compare to the classes considered
in Theorems \ref{2class}, \ref{2classex} and \ref{discad} constructed
from statistical separation conditions. The following result shows that
the conditions in the present paper are strictly weaker than those in
\cite{GN10a,KNP10} for the case of two fixed H\"{o}lder classes, and
also gives a more statistical explanation of why adaptation is possible
over the classes from (\ref{barsig}).
\begin{proposition} \label{gnkill} Let $t>s$.

\begin{longlist}[(a)]
\item[(a)]
Suppose $f \in\bar\Sigma(\varepsilon, s, \ell_n)$ for some fixed
$\varepsilon>0$. Then $\inf_{g \in\Sigma(t)}\|f-g\|_\infty\ge c
2^{-\ell
_n s}$ for some constant $c\equiv c(\varepsilon, B, s,t, K)$.
Moreover, if
$2^{-\ell_ns}/ r_n(s) \to\infty$ as $n \to\infty$, so in particular
in the adaptive case as in (\ref{R}), then, for every $L_0>0$,
\[
\bar\Sigma(\varepsilon, s, \ell_n) \subset\tilde\Sigma(s,L_0 r_n(s))
\]
for $n \ge n_0(\varepsilon, B, s,t, L_0,K)$ large enough.

\item[(b)] If $\ell_n$ is s.t. $2^{-\ell_ns}/ r_n(s) \to\infty$ as $n \to
\infty$, so in particular in the adaptive case (\ref{R}), then
$\forall
L_0'>0, \varepsilon>0$ the set
\[
\tilde\Sigma(s, L_0' r_n(s)) \setminus\bar\Sigma(\varepsilon,s,
\ell_n)
\]
is nonempty for $n \ge n_0(s, t, K, B, L_0')$ large enough.
\end{longlist}
\end{proposition}

%s2.4 ###
\subsection{A Bayesian perspective}

Instead of analyzing the topological capacity of the set removed, one
can try to quantify its size by some measure on the H\"{o}lder space
$\mathcal C^s$. As there is no translation-invariant measure available
we consider certain probability measures on $\mathcal C^s$ that have a
natural interpretation as nonparametric Bayes priors.

Take any $S$-regular wavelet basis $\{\phi_k, \psi_{lk}\dvtx k \in
\mathbb
Z, l \in\mathbb N\}$ of $L^2([0,1]), S>s$. The wavelet
characterization of $\mathcal C^s([0,1])$ motivates to distribute the
basis functions $\psi_{lk}$'s randomly on $\Sigma(s,B)$ as follows: take
$u_{lk}$ i.i.d. uniform random variables on $[-B,B]$ and define the
random wavelet series
\[
U_s(x) = 1+ \sum_{l=J}^\infty\sum_k 2^{-l(s+1/2)} u_{lk} \psi_{lk}(x),
\]
which converges uniformly almost surely. It would be possible to set
$J=0$ and replace $1$ by $\sum_k u_{0k}\phi_k$ below, but to stay
within the density framework we work with this minor simplification,
for which $\int_0^1 U_s(x)\,dx =1$ as well as $U_s \ge0$ almost surely
if $J \equiv J(\|\psi\|_\infty, B,s)$ is chosen large enough. Conclude
that $U_s$ is a random density that satisfies
\[
\|U_s\|_{s, \infty} \le\max\Bigl(1, \sup_{k, l \ge J}|u_{lk}|\Bigr)
\le B\qquad\mbox{a.s.},
\]
so its law is a natural prior on $\Sigma(s,B)$ that uniformly
distributes suitably scaled wavelets on $\Sigma(s)$ around its
expectation $EU_s=1$.
\begin{proposition} \label{unif}
Let $K$ be the wavelet projection kernel associated to a $S$-regular
wavelet basis $\phi, \psi$ of $L^2([0,1])$, $S>s$, and let
$\varepsilon>0,
j \ge0$. Then
\[
\Pr\{\|K_j(U_s)-U_s\|_\infty< \varepsilon B 2^{-js} \} \le
e^{-\log(1/\varepsilon){2^j}}.
\]
\end{proposition}

By virtue of part (a) of Proposition \ref{gnkill} the same bound can be
established, up to constants, for the probability of the sets $\Sigma
(s) \setminus\tilde\Sigma(s, \rho_n)$ under the law of $U_s$.

Similar results (with minor modifications) could be proved if one
replaces the $u_{lk}$'s by i.i.d. Gaussians, which leads to measures
that have a structure similar to Gaussian priors used in Bayesian
nonparametrics; see, for example, \cite{VZ08}. If we choose $j$ at the
natural frequentist rate $2^j \simeq n^{1/(2s+1)}$, then the bound in
Proposition \ref{unif} becomes $e^{-Cn \delta_n^2(s)}, \delta_n(s) =
n^{-s/(2s+1)}$, where $C>0$ can be made as large as desired by choosing
$\varepsilon$ small enough. In view of (2.3) in Theorem 2.1 in
\cite{GGV00} one could therefore heuristically conclude that the exceptional
sets are ``effective null-sets'' from the point of view of Bayesian
nonparametrics.

%s2.5 ###
\subsection{Adaptive confidence bands for collections of H\"{o}lder
classes} \label{cont}

The question arises of how Theorem \ref{2class} can be extended to
adaptation problems over collections of H\"{o}lder classes whose
smoothness degree varies in a fixed interval $[r,R] \subset(0,\infty
)$. A fixed finite number of H\"{o}lder classes can be handled by a
straightforward extension of the proof of Theorem \ref{2class}. Of more
interest is to consider a continuum of smoothness parameters---adaptive
estimators that attain the minimax sup-norm risk over each element of
the collection $\bigcup_{0<s \le R} \Sigma(s)$ exist; see Theorem
\ref{supad} below. Following Theorem \ref{2class} a first approach
might seem to introduce analogues of the sets $\tilde\Sigma(s,
\rho_n)$ as
\[
\Bigl\{f \in\Sigma(s)\dvtx\inf_{g \in\Sigma(t)}\|g-f\|_\infty\ge\rho
_n(s)\mbox{ }\forall t>s \Bigr\}.
\]
However this does not make sense as the sets $\{\Sigma(t)\}_{t>s}$ are
\mbox{$\|\cdot\|_\infty$}-dense in $\Sigma(s)$, so that so-defined $\tilde
\Sigma(s,\rho_n(s))$ would be empty [unless $\rho_n(s)=0$]. Rather one
should note that any adaptation problem with a continuous smoothness
parameter $s$ and convergence rates that are polynomial in $n$ can be
recast as an adaptation problem with a discrete parameter set whose
cardinality grows logarithmically in~$n$. Indeed let us dissect $[r,R]$
into $|\mathcal S_n| \simeq\log n$ points
\[
\mathcal S_n:=\mathcal S_n(\zeta)= \{s_i, i=1,\ldots, \mathcal
|\mathcal S_n|\}
\]
that include $r\equiv s_1,R \equiv s_{|\mathcal S_n|}$,
$s_i<s_{i+1}$ $\forall i$, and each of which has at most $2 \zeta/\log n$
and at least $\zeta/\log n$ distance to the next point, where $\zeta>0$
is a fixed constant. A~simple calculation shows
%
%e2.15 ###
%
\begin{equation} \label{dad}
r_n(s_i) \le C r_n(s)
\end{equation}
for some constant $C=C(\zeta,R)$ and every $s_i \le s < s_{i+1}$, so
that any estimator that is adaptive over $\Sigma(s), s \in\mathcal
S_n$, is also adaptive over $\Sigma(s), s \in[r,R]$.

After this discretization we can define
\[
\tilde\Sigma(s, \rho_n(s), \mathcal S_n) = \Bigl\{f \in\Sigma(s)\dvtx
\inf
_{g \in\Sigma(t)}\|g-f\|_\infty\ge\rho_n(s)\mbox{ }\forall t>s, t
\in
\mathcal S_n \Bigr\},
\]
where $\rho_n(s)$ is a sequence of nonnegative integers. We are
interested in the existence of adaptive confidence bands over
\[
\Sigma(R) \cup\biggl(\bigcup_{s \in\mathcal S_n \setminus\{R\}}
\tilde\Sigma(s, \rho_n(s), \mathcal S_n)\biggr)
\]
under sharp conditions on $\rho_n(s)$.

Let us first address lower bounds, where we consider $T=[0,1]$ for
simplicity. Theorem \ref{2class} cannot be applied directly since the
smoothness index $s$ depends on $n$ in the present setting, and any two
$s,s' \in\mathcal S_n$ could be as close as $\zeta/\log n$ possibly.
If the constant $\zeta$ is taken large enough (but finite) one can
prove the following result.
\begin{theorem}[(Lower bound)]
\label{2classex} Let $T=[0,1], L \ge1$ and $0<\alpha<1/3$ be given, and
let $\mathcal S_n(\zeta)$ be a grid as above. Let $s<s'$ be any two
points in $\mathcal S_n(\zeta)$ and suppose that $C_n$ is a confidence
band that is asymptotically honest with level $\alpha$ over
\[
\Sigma(s') \cup\tilde\Sigma(s, \rho_n(s), \mathcal S_n),
\]
and that is adaptive in the sense that
\[
\sup_{f \in\Sigma(s')}E_f |C_n| \le Lr_n(s'),\qquad \sup_{f \in\tilde
\Sigma(s, \rho_n(s), \mathcal S_n)}E_f|C_n| \le Lr_n(s)
\]
for every $n$ large enough. Then if $\zeta:=\zeta(R,B,L,\alpha)$ is a
large enough but finite constant, we necessarily have
\[
\liminf_n \frac{\rho_n(s)}{r_n(s)}>0.
\]
\end{theorem}

A version of Theorem \ref{2classex} for $T=\mathbb R$ can be proved as
well, by natural modifications of its proof.

To show that adaptive procedures exist if $B$ is known define
\[
\tilde\Sigma_n(s) := \Bigl\{f \in\Sigma(s)\dvtx\inf_{g \in\Sigma
(t)}\|
g-f\|_\infty\ge L_0r_n(s) \mbox{ }\forall t \in\mathcal S_n, t>s
\Bigr\},
\]
where $s$ varies in $[r,R)$, and where $L_0>0$. Setting $\tilde\Sigma
_n(R) \equiv\Sigma(R)$ for notational convenience, we now prove that
an adaptive and honest confidence band exists, for $L_0$ large enough,
over the class
\[
\mathcal P_n(L_0) := \mathcal P(\mathcal S_n,B,L_0,n) := \bigcup_{s
\in
\mathcal S_n} \tilde\Sigma_n(s).
\]
Analyzing the limit set (as $n \to\infty$) of $\mathcal P_n(L_0)$, or
a direct comparison to the continuous scale of classes in (\ref
{barsig}), seems difficult, as $\mathcal S_n$ depends on $n$ now. Note,
however, that one can always choose $\{\mathcal S_n\}_{n \ge1}$ in a
nested way, and $\zeta$ large enough, such that $\mathcal P_n(L_0)$
contains, for every $n$, any fixed finite union (over $s$) of sets of
the form $\bar\Sigma(\varepsilon,s, \ell_n)$ (using Proposition
\ref{gnkill}).
\begin{theorem}[(Existence of adaptive bands)]
\label{discad}
Let $X_1,\ldots,X_n$ be i.i.d. random variables on $T=[0,1]$ or
$T=\mathbb
R$ with density $f \in\mathcal P_n(L_0)$ and suppose $B,r,R, 0<\alpha
<1$ are given. Then, if $L_0$ is large enough depending only on $B$, a
confidence band $C_n=C_n(B,r,R,\alpha; X_1,\ldots, X_n)$ can be
constructed such that
\[
\liminf_n \inf_{f \in\mathcal P_n(L_0)} {\Pr}_f (f \in C_n) \ge
1-\alpha
\]
and, for every $s \in\mathcal S_n, n \in\mathbb N$ and some constant
$L'$ independent of $n$,
%
%e2.16 ###
%
\begin{equation} \label{adc}
\sup_{f \in\tilde\Sigma_n(s)}E_f|C_n| \le L' r_n(s).
\end{equation}
\end{theorem}

%s3 ###
\section{Proofs of remaining results} \label{rem}

\mbox{}

\begin{pf*}{Proof of Proposition \ref{nob}}
On the events $\{f \in C_n\}$ we can find a random density $T_n \in
C_n$ depending only on $C_n$ such that $\{|C_n| \le D, f \in C_n\}
\subseteq\{\|T_n-f\|_\infty\le D \}$ for any $D>0$, and negating this
inclusion we have
\[
\{|C_n| > D\} \cup\{f \notin C_n\} \supseteq\{\|T_n-f\|_\infty> D \}
\]
so that ${\Pr}_f (|C_n| >D) \ge{\Pr}_f (\|T_n-f\|_\infty> D) - {\Pr
}_f(f \notin C_n)$. Thus, using coverage of the band
\begin{eqnarray*}
&& \liminf_n \sup_{f \in\Sigma(r,B)}{\Pr}_f \bigl(|C_n| > cB^p r_n(r)
\bigr) \\
&&\qquad \ge\liminf_n \sup_{f \in\Sigma(r,B)}{\Pr}_f \bigl(\|T_n-f\|_\infty
> cB^p r_n(r) \bigr) - \alpha.
\end{eqnarray*}
The limit inferior in the last line equals $1$ as soon as $c>0$ is
chosen small enough depending only on $r,p$ in view of Theorem 1 in
\cite{KN99}; see also page 1114 as well as Lemma A.2 in that paper.
Taking $\liminf$'s in the inequality
\[
\sup_{f \in\Sigma(r,B)} \frac{E_f|C_n|}{r_n(r)} \ge cB^p \sup_{f
\in
\Sigma(r,B)}{\Pr}_f \bigl(|C_n| > cB^p r_n(r) \bigr)
\]
gives the result.
\end{pf*}
\begin{pf*}{Proof of Proposition \ref{gnkill}}
(a) Observe first that for every $l_0 \ge\ell_n$,
\[
\|\psi\|_\infty\sum_{l \ge l_0} 2^{l/2} \sup_{k}|\beta_{lk}(f)|
\ge\|
K_{l_0}(f)-f\|_\infty\ge\varepsilon2^{-l_0s}.
\]
Let\vspace*{1pt} $N$ be a fixed integer, and let $\ell_n' \ge\ell_n$ be a sequence
of integers to be chosen later. Then for some $\bar l \in[\ell'_n,
\ell'_n+N-1]$
\begin{eqnarray*}
\sup_k |\beta_{\bar l k}(f)| &\ge& \frac{1}{N} \sum_{l= \ell
'_n}^{\ell
'_n+N-1} \sup_k |\beta_{lk}(f)| \\
&\ge& \frac{2^{-(\ell'_n+N)/2}}{N} \Biggl(\sum_{l=\ell'_n}^\infty
2^{l/2} \sup_k |\beta_{lk}(f)| - \sum_{l=\ell'_n+N}^\infty2^{l/2}
\sup
_k |\beta_{lk}(f)| \Biggr) \\
&\ge& \frac{2^{-(\ell'_n+N)/2}}{N} \biggl(\frac{\varepsilon}{\|\psi\|
_\infty
} 2^{-\ell'_ns} - c(B,s)2^{-(\ell'_n+N)s} \biggr) \\
& \ge& \frac{2^{-(\ell'_n+N)/2}}{2 \|\psi\|_\infty N} \varepsilon
2^{-\ell
'_n s} \ge d(\varepsilon,B, \psi,s) 2^{-\ell'_n(s+1/2)}
\end{eqnarray*}
for some $d(\varepsilon, B, \psi,s)>0$ if $N$ is chosen large enough but
finite depending only on $\varepsilon,B, \psi,s$. From (\ref{lbdu}) we
thus have, for any $t>s$,
\begin{eqnarray*}
\inf_{g \in\Sigma(t)}\|f-g\|_\infty&\ge& \inf_{g \in\Sigma(t)}
\sup
_{l \ge\ell'_n, k} 2^{l/2} |\beta_{lk}(f-g)| \\
&\ge& d(\varepsilon,B, \psi,s) 2^{-\ell'_n s} - \sup_{g \in\Sigma(t)}
\sup_{l \ge\ell'_n,k} 2^{l/2}|\beta_{lk}(g)| \\
&\ge& d(\varepsilon,B, \psi,s) 2^{-\ell'_n s} - B2^{-\ell'_n t}\\
&\ge&
c(\varepsilon,B,s,t, \psi) 2^{-\ell_n s},
\end{eqnarray*}
where we have chosen $\ell_n'$ large enough depending only on $B,s,t,
d(\varepsilon, B, \psi,s)$ but still of order $O(\ell_n)$. This completes
the proof of the first claim. The second claim is immediate in view of
the definitions.

(b) Take $f = f_0 + 2^{-\ell_n(s+1/2)} \psi_{\ell_n m}$ for some $m$.
Then $\|f\|_{s, \infty} \le1$ so $f \in\Sigma(s,B)$ and the estimate
in the last display of the proof of part (a) implies
\[
\inf_{g \in\Sigma(t)}\|f-g\|_\infty\ge c 2^{-\ell_ns} \ge L_0' r_n(s)
\]
for $n$ large enough depending only on $B,s,t, L_0', \psi$. On the
other hand\break $\|K_{\ell_n+1}(f)-f\|_\infty=0$ so $f \notin\bar\Sigma
(\varepsilon, s, \ell_n)$ for any $\varepsilon>0$.
\end{pf*}
\begin{pf*}{Proof of Proposition \ref{unif}}
Using (\ref{lbdu}) we have
\[
\|K_j(U_s)-U_s\|_\infty\ge\|\psi\|_1^{-1} \sup_{l \ge j,
k}2^{l/2}|\beta_{lk}(U_s)| \ge\|\psi\|_1^{-1} 2^{-js} \max
_{k=1,\ldots,
2^j}|u_{jk}|.
\]
The variables $u_{jk}/B$ are i.i.d. $U(-1,1)$ and so the $U_k$'s, $U_k
:= |u_{jk}/B|$, are i.i.d. $U(0,1)$ with maximum equal to the largest
order statistic $U_{(2^j)}$. Deduce
\[
\Pr\bigl(\|K_j(U_s)-U_s\|_\infty< \varepsilon B 2^{-js} \bigr) \le\Pr
\bigl(U_{(2^j)} < \varepsilon\bigr) = \varepsilon^{2^j}
\]
to complete the proof.
\end{pf*}
\begin{pf*}{Proof of Theorem \ref{2classex}}
The proof is a modification of the ``necessity part'' of Theorem \ref
{2class}. Let us assume w.l.o.g. $B \ge2, R \ge1$, let us write, in
slight abuse of notation, $s_n, s'_n$ for $s, s'$ throughout this proof
to highlight the dependence on $n$ and choose $j_n(s_n) \in\mathbb N$
such that
\[
(n/ \log n)^{1/(2R+1)} \le c_0 (n/ \log n)^{1/(2s_n+1)} \le
2^{j_n(s_n)} \le(n/ \log n)^{1/(2s_n+1)}
\]
holds for some $c_0>1/(2R+1)^{1/(2R+1)}$ and every $n$ large enough. We
shall assume that $\zeta$ is any fixed number satisfying
\[
\zeta> (4R+2) \max\biggl(\log_2 \bigl((4R+2)B\bigr),(2R+1) \log\frac{(4R+2)
L}{\alpha}\biggr)
\]
in the rest of the proof, and we shall establish $\liminf_n (\rho
_n(s_n)/Lr_n(s_n^+))>0$, where $s_n^+>s_n$ is the larger ``neighbor'' of
$s_n$ in $\mathcal S_n$. This completes the proof since $\liminf_n
r_n(s_n^+)/r_n(s_n) \ge c(\zeta)>0$ by definition of the grid.

Assume thus by way of contradiction that $\liminf_n (\rho
_n(s_n)/Lr_n(s_n^+))=0$ so that, by passing to a subsequence of $n$ if
necessary, $\rho_n(s_n) \le Lr_n(s_n^+)+\delta$ for every $\delta>0$
and every $n=n(\delta)$ large enough. Let $\varepsilon:= 1/(2R+1)$ and
define
\[
f_0=1,\qquad f_m= f_0 + \varepsilon2^{-j(s_n+1/2)} \psi_{jm},\qquad
m=1,\ldots, M,
\]
as in the proof of Theorem \ref{2class}, $c_0'2^j \le M \le2^{j}$,
$c_0'>0$. Then $f_m \in\Sigma(s_n)$ for every $j \ge j_0$ where $j_0$
can be taken to depend only on $r,R, B, \psi$. Moreover for $j \ge
(\log n)/(4R+2)$ we have, using (\ref{lbdu}) and the assumption on
$\zeta$, for any $g \in\Sigma(t), t \in\mathcal S_n, t>s_n$, and
every $m$
%
%e3.1 ###
%
\begin{eqnarray} \label{agl}
\|f_m-g\|_\infty&\ge& \sup_{l \ge0, k} 2^{l/2} |\beta
_{lk}(f_m)-\beta
_{lk}(g)| \nonumber\\
& \ge& \varepsilon2^{-js_n} - 2^{j/2}|\beta_{jk}(g)| \ge\varepsilon
2^{-js_n} - B2^{-jt} \\
&\ge& 2^{-js_n}(\varepsilon-B2^{-j\zeta/\log n})
\ge\frac{\varepsilon}{2}2^{-js_n}.\nonumber
\end{eqnarray}
We thus see that
\[
f_m \in\tilde\Sigma\biggl(s_n, \frac{\varepsilon}{2}2^{-js_n}, \mathcal
S_n\biggr) \qquad\forall m=1,\ldots,M,
\]
for every $j \ge J_0 := \max(j_0, (\log n)/(4R+2))$. Take now $j
\equiv
j_n(s_n)$ which exceeds $J_0$ for $n$ large enough, and conclude
%
%e3.2 ###
%
\begin{eqnarray}\label{estlow}
\frac{\varepsilon}{2} 2^{-j_n(s_n)s_n} &\ge& \frac{\varepsilon}{2} r_n(s_n)
\ge\frac{\varepsilon}{2} \frac{\alpha}{L} \frac{L}{\alpha}
(e^{\zeta
/2})^{{1}/{(2R+1)^2}} r_n(s_n^+) \nonumber\\[-8pt]\\[-8pt]
&\ge& \frac{L}{\alpha} r_n(s_n^+) \ge\rho_n(s_n)\nonumber
\end{eqnarray}
for $n$ large enough, where we have used the definition of the grid
$\mathcal S_n$, of $\varepsilon$, the assumption on $\zeta$ and the
hypothesis on $\rho_n$. Summarizing $f_m \in\tilde\Sigma(s_n, \rho
_n(s_n), \mathcal S_n)$ for every $m=1,\ldots,M$ and every $n \ge n_0,
n_0=n_0(r, R, B, \psi)$.

Suppose now $C_n$ is a confidence band that is adaptive and
asymptotically honest over $\Sigma(s'_n) \cup\tilde\Sigma(s_n, \rho
_n(s_n), \mathcal S_n)$, and consider testing $H_0\dvtx f = f_0$ against
$H_1\dvtx f \in\{f_1,\ldots, f_M\}=:\mathcal M$. Define a test $\Psi
_n$ as
follows: if no $f_m \in C_n$ then $\Psi_n=0$, but as soon as one of the
$f_m$'s is contained in $C_n$ then $\Psi_n=1$. Now since $r_n(s'_n)
\le
r_n(s^+_n)$ and using (\ref{agl}), (\ref{estlow}) we have
\begin{eqnarray*}
{\Pr}_{f_0}(\Psi_n \ne0) &=& {\Pr}_{f_0} (f_m \in C_n \mbox{ for some
} m) \\
&\le& {\Pr}_{f_0} (\|f_m-f_0\|_\infty\le|C_n| \mbox{ for some } m) +
\alpha+ o(1) \\
&\le& {\Pr}_{f_0}\bigl(|C_n| \ge(L/\alpha)r_n(s_n^+)\bigr) + \alpha+o(1)
\\
&\le& \alpha r_n(s'_n)/r_n(s^+_n) + \alpha+o(1) \le2 \alpha+ o(1).
\end{eqnarray*}
Under any alternative $f_m \in\tilde\Sigma(s_n)$, invoking honesty of
the band we have
\[
P_{f_m}(\Psi_n = 0) = {\Pr}_{f_m}(\mbox{no } f_k \in C_n) \le{\Pr
}_{f_m} (f_m \notin C_n) \le\alpha+ o(1)
\]
so that summarizing we have
\[
\limsup_n\Bigl(E_{f_0}\Psi_n + \sup_{f \in\mathcal M}E_{f}(1-\Psi
_n)\Bigr) \le3 \alpha<1.
\]
But\vspace*{1pt} this has led to a contradiction by the same arguments as in the
proof\break of Theorem \ref{2class}, noting in the last step that $n\gamma
^2_j = \varepsilon^2 n2^{-j_n(s_n)(2s_n+1)} \le(\varepsilon^2/\break(c_0)^{2R+1})
\log n$ and thus
\[
\frac{e^{n\gamma_j^2}-1}{M} \le\frac{1}{c_0' c_0} e^{
({\varepsilon
^2}/{(c_0)^{2R+1}}) \log n} \biggl(\frac{\log n}{n}\biggr)^{1/(2R+1)} = o(1)
\]
since $1/(2R+1)=\varepsilon<c_0^{2R+1}$.
\end{pf*}
\begin{pf*}{Proof of Theorem \ref{discad}} We shall only prove the more
difficult case \mbox{$T=\mathbb R$}. Let $j_i$ be such that $2^{j_i} \simeq
(n/\log n)^{1/(2s_i+1)}$, let $f_n(j)$ be as in (\ref{est}) based on
wavelets of regularity $S>R$ and define test statistics
\[
\hat d_n(i):=\inf_{g \in\Sigma(s_{i+1})}\|f_n(j_i)-g\|_\infty,\qquad
i=1,\ldots, |\mathcal S_n|-1.
\]
Recall further $\sigma(j)$ from (\ref{sigma}) and, for a constant $L$
to be chosen below, define tests
\[
\Psi(i)=
\cases{0, &\quad if $\hat d_n(i) \le L \sigma(j_i)$, \cr
1, &\quad otherwise,}
\]
to accept $H_0\dvtx f \in\Sigma(s_{i+1})$ against the alternative
$H_1\dvtx f
\in\tilde\Sigma_n(s_i)$. Starting from the largest model we first
test $H_0\dvtx f \in\Sigma(s_2)$ against $H_1\dvtx f \in\tilde
\Sigma_n(r)$.
If $H_0$ is rejected we set $\hat s_n = r$, otherwise we proceed to
test $H_0\dvtx f \in\Sigma(s_3)$ against $H_1\dvtx f \in\tilde
\Sigma
_n(s_2)$ and iterating this procedure downwards we define $\hat s_n$ to
be the first element $s_i$ in $\mathcal S$ for which $\Psi(i)=1$
rejects. If no rejection occurs set $\hat s_n=R$.

For $f \in\mathcal P_n(L_0)$ define $s_{i_0}:=s_{i_0}(f)= \max\{s \in
\mathcal S_n\dvtx f \in\tilde\Sigma_n(s) \}$.
\begin{lemma}
We can choose the constants $L$ and then $L_0$ depending only on $B,
\phi, \psi$ such that
\[
\sup_{f \in\mathcal P_n(L_0)}{\Pr}_f\bigl(\hat s_n \ne s_{i_0}(f)\bigr) \le Cn^{-2}
\]
for some constant $C$ and every $n$ large enough.
\end{lemma}
\begin{pf}
If $\hat s_n < s_{i_0}$, then the test $\Psi(i)$ has rejected for some
$i < i_0$. In this case $f \in\tilde\Sigma_n(s_{i_0}) \subset\Sigma
(s_{i_0}) \subseteq\Sigma(s_{i+1})$ for every $i<i_0$, and thus,
proceeding as in (\ref{bias}) and using Proposition \ref{ineq} below,
we have for $L$ and then $d$ large enough depending only on $B,K$
\begin{eqnarray*}
{\Pr}_f(\hat s_n< s_{i_0}) &=& {\Pr}_f\biggl(\bigcup_{i < i_0}\Bigl\{\inf
_{g \in\Sigma(s_{i+1})} \|f_n(j_i)-g\|_\infty>L \sigma(j_i)\Bigr\}
\biggr) \\
&\le& \sum_{i<i_0} {\Pr}_f\bigl(\|f_n(j_i)-E_ff_n(j_i)\|_\infty>(L-b)
\sigma(j_i)\bigr) \\
&\le& C' |\mathcal S_n| e^{-d \log n} \le C n^{-2}.
\end{eqnarray*}

On the other hand if $\hat s_n > s_{i_0}$ (ignoring the trivial case
$s_{i_0} =R$), then $\Psi(i_0)$ has accepted despite $f \in\tilde
\Sigma_n(s_{i_0})$. Thus, using $r_n(s_{i_0}) \ge c \sigma(j_{i_0})$
for some $c=c(B)$ and proceeding as in (\ref{bias}) we can bound ${\Pr
}_f(\hat s_n > s_{i_0})$ by
\begin{eqnarray*}
&& {\Pr}_f \Bigl(\inf_{g \in\Sigma(s_{i_0+1})} \|f_n(j_{i_0})-g\|
_\infty\le L \sigma(j_{i_0}) \Bigr) \\
&&\qquad\le{\Pr}_f \Bigl(\inf_{g \in\Sigma(s_{i_0+1})} \|f-g\|_\infty- \|
f_n(j_{i_0})-E_ff_n(j_{i_0})\|_\infty\\
&&\qquad\quad\hspace*{87pt}{} - \|E_ff_n(j_{i_0})-f\|_\infty\le L \sigma(j_{i_0})
\Bigr) \\
&&\qquad \le{\Pr}_f \bigl(L_0r_n(s_{i_0}) - \|K_{j_{i_0}}(f)-f\|_\infty- L
\sigma(j_{i_0}) \le\|f_n(j_{i_0})-E_ff_n(j_{i_0}) \|_\infty\bigr) \\
&&\qquad \le{\Pr}_f \bigl(\|f_n(j_{i_0})-E_ff_n(j_{i_0}) \|_\infty\ge
(cL_0-L-b) \sigma(j_{i_0}) \bigr) \\
&&\qquad \le c'e^{-c'j_{i_0}} \le C/n^2
\end{eqnarray*}
for $L_0$ and then also $c'>0$ large enough, using Proposition \ref
{ineq} below.
\end{pf}

Take now $\hat f_n$ to be an estimator of $f$ that is adaptive in
sup-norm loss over $\bigcup_{s \in[r,R]} \Sigma(s)$ as in Theorem
\ref
{supad} below and define the confidence band
\[
C_n = \hat f_n \pm M \biggl(\frac{\log n}{n} \biggr)^{{\hat
s_n}/({2\hat s_n +1})},
\]
where $M$ is chosen below. For $f \in\tilde\Sigma_n(s_{i_0})$ the
lemma implies
\begin{eqnarray*}
E_f|C_n| &\le& 2M \biggl(\frac{\log n}{n} \biggr)^{
{s_{i_0}}/({2s_{i_0} +1})} + 2M \biggl(\frac{\log n}{n} \biggr)^{
{r}/({2r +1})} \times{\Pr}_f (\hat s_n <s_{i_0}) \\
&\le& C(M) \biggl(\frac{\log n}{n} \biggr)^{{s_{i_0}}/({2s_{i_0} +1})},
\end{eqnarray*}
so this band is adaptive.

For coverage, we have, again from the lemma and Markov's inequality
\begin{eqnarray*}
{\Pr}_f (f \in C_n ) &=& {\Pr}_f \bigl(\|\hat f_n -f\|
_\infty\le Mr_n(\hat s_n) \bigr) \\
&\ge& 1-{\Pr}_f \bigl(\|\hat f_n -f\|_\infty> Mr_n(s_{i_0}) \bigr) -
\Pr(\hat s_n >s_{i_0} ) \\
&\ge& 1- \frac{E_f\|\hat f_n -f\|_\infty}{Mr_n(s_{i_0})} - \frac
{C}{n^2} \\
&\ge&1 - \frac{D(B, R,r)}{M} - \frac{C}{n^2},
\end{eqnarray*}
which is greater than or equal to $1-\alpha$ for $M$ and $n$ large
enough depending only on $B,R,r$.
\end{pf*}
\begin{pf*}{Proof of part \textup{(c)} of Theorem \ref{2class}}
The analog of case (b) is immediate. The analog of part (a) requires
the following modifications: set again $f_0=1$ on $[0,1]$, $0\le j_n' <
j_n$ to be chosen below, and define
\[
f_m:= f_0 + B 2^{-j_n'(s+1/2)} \psi_{j_n'm_0} + \varepsilon
2^{-j_n(r+1/2)} \psi_{j_nm},
\]
where $m=1,\ldots, M \simeq2^j$, all $\psi_{lk}$'s are Daubechies wavelets
supported in the interior of $[0,1]$ and where $m_0\ne m$ is chosen
such that $\psi_{j_n'm_0}$ and $\psi_{j_nm}$ have disjoint support for
every $m$ (which is possible for $j_n, j_n'$ large enough since
Daubechies wavelets have localized support). Recalling $j_n^*$ from the
proof of part (a), we can choose $j_n', j_n$ in such a way that
$j'_n<j_n, 2^{-j_nr} =o(2^{-j_n^*r})$,
\[
f_m \in\tilde\Sigma(r, \rho_n)\qquad\forall m,\qquad f'_0 := f_0+B
2^{-j_n'(s+1/2)} \psi_{j_n'm_0} \in\tilde\Sigma\bigl(s, (B/2) r_n(s)\bigr)
\]
for every $n \ge n_0, n_0=n_0(s,r,B,\varepsilon, \psi)$. Now if $C_n$
is a
confidence band that is adaptive and honest over $\tilde\Sigma(s,
r_n(s)) \cup\tilde\Sigma(r, \rho_n)$ consider testing $H_0\dvtx f = f'_0
$ against $H_1\dvtx f \in\{f_1,\ldots, f_M\}=:\mathcal M$. The same
arguments as before (\ref{bul}) show that there exists a test $\Psi_n$
such that $\limsup_n(E_{f_0}\Psi_n + \sup_{f \in\mathcal
M}E_{f}(1-\Psi
_n)) \le2 \alpha<1$ along a subsequence of $n$, a claim that leads to
a contradiction since we can lower bound the error probabilities of any
test as in the original proof above, the only modification arising in
the bound for the likelihood ratio. Let $P_0'$ be the $n$-fold product
probability measure induced by the density $f_0'$ and set $Z= (1/M)
\sum
_{m=1}^M (dP_m/dP'_0)$. We suppress now the dependence of $j_n$ on $n$
for notational simplicity, and define shorthand $\gamma_j =
\varepsilon
2^{-j(r+1/2)}$, $\kappa_j = B 2^{-j'(s+1/2)}$. To bound
$E_{f'_0}(Z-1)^2$ we note that, using orthonormality of the $\psi
_{jm}$'s, that $\int\psi_{jm}=0$ and that $\psi_{j'm_0}$ has disjoint
support with $\psi_{jm}, m=1,\ldots, M$, we have ($m\ne m'$)
\begin{eqnarray*}
\int\frac{\psi_{lm}\psi_{lm'}}{(1+\kappa_j \psi_{j'm_0})^2}f_0' &=&
\int
\psi_{jm}\psi_{jm'}=0,
\\
\int\frac{\psi_{jm}}{1+\kappa_j \psi_{j'm_0}}f_0' &=& \int\psi
_{jm} =0,\\
\int\frac{\psi^2_{jm}}{(1+\kappa_j \psi_{j'm_0})^2}f_0'&=&\int\psi
^2_{jm}=1.
\end{eqnarray*}
The identities in the last display can be used to bound
$E_{f'_0}(Z-1)^2$ by
\begin{eqnarray*}
&&\frac{1}{M^2} \int_{[0,1]^n} \Biggl(\sum_{m=1}^M \Biggl(\prod_{i=1}^n \frac
{f_m(x_i)}{f_0'(x_i)} -1\Biggr) \Biggr)^2 \prod_{i=1}^n f_0'(x_i) \,dx \\
&&\qquad= \frac{1}{M^2} \sum_{m=1}^M \int_{[0,1]^n} \Biggl( \prod_{i=1}^n
\biggl(1+\frac{\gamma_j \psi_{jm}(x_i)}{1+\kappa_j \psi_{j'm_0}(x_i)}\biggr)
-1 \Biggr)^2 \prod_{i=1}^nf'_0(x_i)\,dx \\
&&\qquad= \frac{1}{M^2} \sum_{m=1}^M \biggl( \biggl( \int_{[0,1]} \biggl(1+\frac
{\gamma_j \psi_{jm}(x_i)}{1+\kappa_j \psi_{j'm_0}(x_i)}\biggr)^2
f_0'(x)\,dx \biggr)^n -1 \biggr) \\
&&\qquad= \frac{1}{M} \bigl((1+\gamma_j^2)^n -1 \bigr) \le\frac{e^{n \gamma
_j^2}-1}{M}.
\end{eqnarray*}
The rest of the proof is as in part (a) of Theorem \ref{2class}.
\end{pf*}

%s3.1 ###
\subsection{Auxiliary results}

The following theorem is due to \cite{GN09a,GN09b,GN10b}. We state a
version that follows from Theorem 4 in \cite{LN10} for $T=\mathbb R$.
In case $T=[0,1]$ it follows from the same proofs. The restriction that
$B$ be known is not necessary but suffices for our present purposes.
\begin{theorem}\label{supad} Let $X_1,\ldots,X_n$ be i.i.d. with
uniformly continuous density $f$ on $T=[0,1]$ or $T=\mathbb R$. Then
for every $r,R, 0<r \le R$ there exists an estimator $\hat f_n(x):=\hat
f_n(x,X_1,\ldots, X_n,B,R)$ such that, for every $s$, $r \le s\le R$,
some constant $D(B,r,R)$ and every $n \ge2$ we have $\sup_{f \in
\Sigma
(s,B,T)}E \|\hat f_n-f\|_\infty\le D(B,r,R) r_n(s)$.
\end{theorem}

The following inequality was proved in \cite{GN09b} (see also page 1167
in \cite{GN10a}) for $T=\mathbb R$ (the case $T=[0,1]$ is similar, in
fact simpler).
\begin{proposition} \label{ineq}
Let $\phi, \psi$ be a compactly supported scaling and wavelet function,
respectively, both $S$-H\"{o}lder for some $S>0$. Suppose $P$ has a
bounded density $f$ and let $f_n(x,j)$ be the estimator from (\ref
{est}). Given $C, C'>0$, there exist finite positive constants $C_1 =
C_1(C,K)$ and $C_2 = C_2(C,C',K)$ such that, if $(n/2^j j) \geq C$ and
$C_1\sqrt{(\|f\|_\infty\vee1) (2^jj/n)} \leq t \leq C'$, then, for
every $n \in\mathbb N$,
\[
{\Pr}_f\Bigl\{\sup_{x \in\mathbb R} |f_n(x,j)-Ef_n(x,j)| \ge t
\Bigr\} \le C_2\exp\biggl(-\frac{nt^2}{C_2(\|f\|_\infty\vee1)2^j}\biggr).
\]
\end{proposition}

\section*{Acknowledgments}
We would like to thank Tony Cai, the Associate Editor, as well as two
anonymous referees for valuable remarks and criticism. We are further
grateful to Adam Bull and Jakob S\"{o}hl for pointing out mistakes in
the first version of this article. We are also indebted to Yannick
Baraud, Lutz D\"{u}mbgen, Vladimir Koltchinskii, Oleg Lepski, Axel
Munk, Dominique Picard, Benedikt P\"{o}tscher, Markus Reiss, Vladimir
Spokoiny and Aad van der Vaart for stimulating discussions on the
subject over the past years. The second author would like to thank the
hospitality of the Caf\'{e}s Br\"{a}unerhof and Florianihof in Vienna
for their hospitality.

%suskaldyti doi

% imsref loaded by lrinkeviciute, 2011-09-28 13:13:59
% imsref loaded by lrinkeviciute, 2011-09-28 13:15:00
%

\printaddresses

\end{document}